\documentclass{amsart}
\usepackage{amsfonts}

\newtheorem{theorem}{Theorem}
\theoremstyle{plain}

\newtheorem{corollary}{Corollary}

\newtheorem{proposition}{Proposition}
\newtheorem{remark}{Remark}

\numberwithin{equation}{section}
\input{tcilatex}

\begin{document}
\title[Bessel and Gr\"{u}ss Inequalities]{Some New Results Related to Bessel
and Gr\"{u}ss Inequalities for Orthogonal Families in Inner Product Spaces}
\author{S.S. Dragomir}
\address{School of Computer Science and Mathematics\\
Victoria University of Technology\\
PO Box 14428, MCMC 8001\\
Victoria, Australia.}
\date{May 19, 2003}
\email{sever.dragomir@vu.edu.au}
\urladdr{http://rgmia.vu.edu.au/SSDragomirWeb.html}
\subjclass[2000]{26D15, 46C05.}
\keywords{Bessel's inequality, Gr\"{u}ss' inequality, Inner product,
Lebesgue integral.}

\begin{abstract}
Some new counterparts of Bessel's inequality for orthornormal families in
real or complex inner product spaces are pointed out. Applications for some
Gr\"{u}ss type inequalities are also empahsized.
\end{abstract}

\maketitle

\section{Introduction\label{s1}}

In \cite{SSDa}, the author has proved the following result which provides
both a Gr\"{u}ss type inequality for orthogonal families of vectors in real
or complex inner products as well as, for $x=y,$ a counterpart of Bessel's
inequality.

\begin{theorem}
\label{t1.1}Let $\left\{ e_{i}\right\} _{i\in I}$ be a family of
orthornormal vectors in $H,$ i.e., $\left\langle e_{i},e_{j}\right\rangle =0$
if $i\neq j$ and $\left\Vert e_{i}\right\Vert =1,$ $i,j\in I,$ $F$ a finite
part of $I,$ $\phi _{i},\gamma _{i},\Phi _{i},\Gamma _{i}\in \mathbb{R}$ $%
\left( i\in F\right) $, and $x,y\in H.$ If either%
\begin{align}
\func{Re}\left\langle \sum_{i=1}^{n}\Phi _{i}e_{i}-x,x-\sum_{i=1}^{n}\phi
_{i}e_{i}\right\rangle & \geq 0,\   \label{1.1} \\
\func{Re}\left\langle \sum_{i=1}^{n}\Gamma
_{i}e_{i}-y,y-\sum_{i=1}^{n}\gamma _{i}e_{i}\right\rangle & \geq 0,  \notag
\end{align}%
or, equivalently,%
\begin{align}
\left\Vert x-\sum_{i\in F}\frac{\Phi _{i}+\phi _{i}}{2}e_{i}\right\Vert &
\leq \frac{1}{2}\left( \sum_{i\in F}\left\vert \Phi _{i}-\phi
_{i}\right\vert ^{2}\right) ^{\frac{1}{2}},  \label{1.2} \\
\left\Vert y-\sum_{i\in F}\frac{\Gamma _{i}+\gamma _{i}}{2}e_{i}\right\Vert
& \leq \frac{1}{2}\left( \sum_{i\in F}\left\vert \Gamma _{i}-\gamma
_{i}\right\vert ^{2}\right) ^{\frac{1}{2}},  \notag
\end{align}%
hold, then we have the inequality%
\begin{align}
0& \leq \left\vert \left\langle x,y\right\rangle -\sum_{i\in F}\left\langle
x,e_{i}\right\rangle \left\langle e_{i},y\right\rangle \right\vert
\label{1.3} \\
& \leq \frac{1}{4}\left( \sum_{i\in F}\left\vert \Phi _{i}-\phi
_{i}\right\vert ^{2}\right) ^{\frac{1}{2}}\cdot \left( \sum_{i\in
F}\left\vert \Gamma _{i}-\gamma _{i}\right\vert ^{2}\right) ^{\frac{1}{2}} 
\notag \\
& \ \ \ \ \ \ \ \ \ \ \ -\left[ \func{Re}\left\langle \sum_{i\in F}\Phi
_{i}e_{i}-x,x-\sum_{i\in F}\phi _{i}e_{i}\right\rangle \right] ^{\frac{1}{2}}
\notag \\
& \ \ \ \ \ \ \ \ \ \ \ \ \times \left[ \func{Re}\left\langle \sum_{i\in
F}\Gamma _{i}e_{i}-y,y-\sum_{i\in F}\gamma _{i}e_{i}\right\rangle \right] ^{%
\frac{1}{2}}  \notag \\
& \leq \frac{1}{4}\left( \sum_{i\in F}\left\vert \Phi _{i}-\phi
_{i}\right\vert ^{2}\right) ^{\frac{1}{2}}\cdot \left( \sum_{i\in
F}\left\vert \Gamma _{i}-\gamma _{i}\right\vert ^{2}\right) ^{\frac{1}{2}}. 
\notag
\end{align}%
The constant $\frac{1}{4}$ is best possible in the sense that it cannot be
replaced by a smaller constant.
\end{theorem}

In the follow up paper \cite{SSDb}, and by the use of a different technique,
the author has pointed out the following result as well:

\begin{theorem}
\label{t2}Let $\left\{ e_{i}\right\} _{i\in I}$, $F,$ $\phi _{i},\gamma
_{i},\Phi _{i},\Gamma _{i}$ and $x,y$ be as in Theorem \ref{t1.1}. If either
(\ref{1.1}) or (\ref{1.2}) holds, then we have the inequality%
\begin{align}
0& \leq \left\vert \left\langle x,y\right\rangle -\sum_{i=1}^{n}\left\langle
x,e_{i}\right\rangle \left\langle e_{i},y\right\rangle \right\vert
\label{1.4} \\
& \leq \frac{1}{4}\left( \sum_{i=1}^{n}\left\vert \Phi _{i}-\phi
_{i}\right\vert ^{2}\right) ^{\frac{1}{2}}\cdot \left(
\sum_{i=1}^{n}\left\vert \Gamma _{i}-\gamma _{i}\right\vert ^{2}\right) ^{%
\frac{1}{2}}  \notag \\
& \ \ \ \ \ \ \ \ \ \ \ \ \ \ \ \ \ -\sum_{i\in F}\left\vert \frac{\Phi
_{i}+\phi _{i}}{2}-\left\langle x,e_{i}\right\rangle \right\vert \left\vert 
\frac{\Gamma _{i}+\gamma _{i}}{2}-\left\langle y,e_{i}\right\rangle
\right\vert  \notag \\
& \leq \frac{1}{4}\left( \sum_{i=1}^{n}\left\vert \Phi _{i}-\phi
_{i}\right\vert ^{2}\right) ^{\frac{1}{2}}\cdot \left(
\sum_{i=1}^{n}\left\vert \Gamma _{i}-\gamma _{i}\right\vert ^{2}\right) ^{%
\frac{1}{2}}.  \notag
\end{align}%
The constant $\frac{1}{4}$ is best possible in the sense that it cannot be
replaced by a smaller constant.
\end{theorem}

It has also been shown that the bounds provided by the second inequality in (%
\ref{1.3}) and the second inequality in (\ref{1.4}) cannot be compared in
general.

\section{A New Counterpart of Bessel's Inequality\label{s2}}

The following counterpart of Bessel's inequality holds.

\begin{theorem}
\label{t2.1}Let $\left\{ e_{i}\right\} _{i\in I}$ be a family of
orthornormal vectors in $H,$ $F$ a finite part of $I,$ and $\phi _{i},\Phi
_{i}$ $\left( i\in F\right) ,$ real or complex numbers such that $\sum_{i\in
F}\func{Re}\left( \Phi _{i}\overline{\phi _{i}}\right) >0.$ If $x\in H$ is
such that either

\begin{enumerate}
\item[(i)] $\func{Re}\left\langle \sum_{i\in F}\Phi _{i}e_{i}-x,x-\sum_{i\in
F}\phi _{i}e_{i}\right\rangle \geq 0;$\newline
or, equivalently,

\item[(ii)] $\left\Vert x-\sum_{i\in F}\frac{\phi _{i}+\Phi _{i}}{2}%
e_{i}\right\Vert \leq \frac{1}{2}\left( \sum_{i\in F}\left\vert \Phi
_{i}-\phi _{i}\right\vert ^{2}\right) ^{\frac{1}{2}};$
\end{enumerate}

holds, then one has the inequality 
\begin{equation}
\left\Vert x\right\Vert ^{2}\leq \frac{1}{4}\cdot \frac{\sum_{i\in F}\left(
\left\vert \Phi _{i}\right\vert +\left\vert \phi _{i}\right\vert \right) ^{2}%
}{\sum_{i\in F}\func{Re}\left( \Phi _{i}\overline{\phi _{i}}\right) }%
\sum_{i\in F}\left\vert \left\langle x,e_{i}\right\rangle \right\vert ^{2}.
\label{2.1}
\end{equation}%
The constant $\frac{1}{4}$ is best possible in the sense that it cannot be
replaced by a smaller constant.
\end{theorem}

\begin{proof}
Firstly, we observe that for $y,a,A\in H,$ the following are equivalent%
\begin{equation}
\func{Re}\left\langle A-y,y-a\right\rangle \geq 0  \label{2.2}
\end{equation}%
and%
\begin{equation}
\left\Vert y-\frac{a+A}{2}\right\Vert \leq \frac{1}{2}\left\Vert
A-a\right\Vert .  \label{2.3}
\end{equation}%
Now, for $a=\sum_{i\in F}\phi _{i}e_{i},$ $A=\sum_{i\in F}\Phi _{i}e_{i},$
we have%
\begin{align*}
\left\Vert A-a\right\Vert & =\left\Vert \sum_{i\in F}\left( \Phi _{i}-\phi
_{i}\right) e_{i}\right\Vert =\left[ \left\Vert \sum_{i\in F}\left( \Phi
_{i}-\phi _{i}\right) e_{i}\right\Vert ^{2}\right] ^{\frac{1}{2}} \\
& =\left( \sum_{i\in F}\left\vert \Phi _{i}-\phi _{i}\right\vert
^{2}\left\Vert e_{i}\right\Vert ^{2}\right) ^{\frac{1}{2}}=\left( \sum_{i\in
F}\left\vert \Phi _{i}-\phi _{i}\right\vert ^{2}\right) ^{\frac{1}{2}},
\end{align*}%
giving, for $y=x,$ the desired equivalence.

Now, observe that%
\begin{multline*}
\func{Re}\left\langle \sum_{i\in F}\Phi _{i}e_{i}-x,x-\sum_{i\in F}\phi
_{i}e_{i}\right\rangle \\
=\sum_{i\in F}\func{Re}\left[ \Phi _{i}\overline{\left\langle
x,e_{i}\right\rangle }+\overline{\phi _{i}}\left\langle x,e_{i}\right\rangle %
\right] -\left\Vert x\right\Vert ^{2}-\sum_{i\in F}\func{Re}\left( \Phi _{i}%
\overline{\phi _{i}}\right) ,
\end{multline*}%
giving, from (i), that%
\begin{equation}
\left\Vert x\right\Vert ^{2}+\sum_{i\in F}\func{Re}\left( \Phi _{i}\overline{%
\phi _{i}}\right) \leq \sum_{i\in F}\func{Re}\left[ \Phi _{i}\overline{%
\left\langle x,e_{i}\right\rangle }+\overline{\phi _{i}}\left\langle
x,e_{i}\right\rangle \right] .  \label{2.4}
\end{equation}

On the other hand, by the elementary inequality%
\begin{equation*}
\alpha p^{2}+\frac{1}{\alpha }q^{2}\geq 2pq,\ \ \alpha >0,\ p,q\geq 0;
\end{equation*}%
we deduce%
\begin{equation}
2\left\Vert x\right\Vert \leq \frac{\left\Vert x\right\Vert ^{2}}{\left[
\sum_{i\in F}\func{Re}\left( \Phi _{i}\overline{\phi _{i}}\right) \right] ^{%
\frac{1}{2}}}+\left[ \sum_{i\in F}\func{Re}\left( \Phi _{i}\overline{\phi
_{i}}\right) \right] ^{\frac{1}{2}}.  \label{2.5}
\end{equation}%
Dividing (\ref{2.4}) by $\left[ \sum_{i\in F}\func{Re}\left( \Phi _{i}%
\overline{\phi _{i}}\right) \right] ^{\frac{1}{2}}>0$ and using (\ref{2.5}),
we obtain%
\begin{equation}
\left\Vert x\right\Vert \leq \frac{1}{2}\frac{\sum_{i\in F}\func{Re}\left[
\Phi _{i}\overline{\left\langle x,e_{i}\right\rangle }+\overline{\phi _{i}}%
\left\langle x,e_{i}\right\rangle \right] }{\left[ \sum_{i\in F}\func{Re}%
\left( \Phi _{i}\overline{\phi _{i}}\right) \right] ^{\frac{1}{2}}},
\label{2.6}
\end{equation}%
which is also an interesting inequality in itself.

Using the Cauchy-Buniakowsky-Schwartz inequality for real numbers, we get%
\begin{align}
\sum_{i\in F}\func{Re}\left[ \Phi _{i}\overline{\left\langle
x,e_{i}\right\rangle }+\overline{\phi _{i}}\left\langle x,e_{i}\right\rangle %
\right] & \leq \sum_{i\in F}\left\vert \Phi _{i}\overline{\left\langle
x,e_{i}\right\rangle }+\overline{\phi _{i}}\left\langle x,e_{i}\right\rangle
\right\vert  \label{2.7} \\
& \leq \sum_{i\in F}\left( \left\vert \Phi _{i}\right\vert +\left\vert \phi
_{i}\right\vert \right) \left\vert \left\langle x,e_{i}\right\rangle
\right\vert  \notag \\
& \leq \left[ \sum_{i\in F}\left( \left\vert \Phi _{i}\right\vert
+\left\vert \phi _{i}\right\vert \right) ^{2}\right] ^{\frac{1}{2}}\left[
\sum_{i\in F}\left\vert \left\langle x,e_{i}\right\rangle \right\vert ^{2}%
\right] ^{\frac{1}{2}}.  \notag
\end{align}%
Making use of (\ref{2.6}) and (\ref{2.7}), we deduce the desired result (\ref%
{2.1}).

To prove the sharpness of the constant $\frac{1}{4},$ let us assume that (%
\ref{2.1}) holds with a constant $c>0,$ i.e., 
\begin{equation}
\left\Vert x\right\Vert ^{2}\leq c\cdot \frac{\sum_{i\in F}\left( \left\vert
\Phi _{i}\right\vert +\left\vert \phi _{i}\right\vert \right) ^{2}}{%
\sum_{i\in F}\func{Re}\left( \Phi _{i}\overline{\phi _{i}}\right) }%
\sum_{i\in F}\left\vert \left\langle x,e_{i}\right\rangle \right\vert ^{2},
\label{2.8}
\end{equation}%
provided $x,$ $\phi _{i},\Phi _{i},i\in F$ satisfies (i).

Choose $F=\left\{ 1\right\} ,$ $e_{1}=e,$ $\left\Vert e\right\Vert =1,$ $%
\phi _{i}=m,$ $\Phi _{i}=M$ with $m,M>0,$ then, by (\ref{2.8}) we get%
\begin{equation}
\left\Vert x\right\Vert ^{2}\leq c\frac{\left( M+m\right) ^{2}}{mM}%
\left\vert \left\langle x,e\right\rangle \right\vert ^{2}  \label{2.9}
\end{equation}%
provided%
\begin{equation}
\func{Re}\left\langle Me-x,x-me\right\rangle \geq 0.  \label{2.10}
\end{equation}%
If $x=me,$ then obviously (\ref{2.10}) holds, and by (\ref{2.9}) we get%
\begin{equation*}
m^{2}\leq c\frac{\left( M+m\right) ^{2}}{mM}m^{2}
\end{equation*}%
giving $mM\leq c\left( M+m\right) ^{2}$ for $m,M>0.$ Now, if in this
inequality we choose $m=1-\varepsilon ,$ $M=1+\varepsilon $ $\left(
\varepsilon \in \left( 0,1\right) \right) ,$ then we get $1-\varepsilon
^{2}\leq 4c$ for $\varepsilon \in \left( 0,1\right) ,$ from where we deduce $%
c\geq \frac{1}{4}.$
\end{proof}

\begin{remark}
\label{r2.2}By the use of (\ref{2.6}), the second inequality in (\ref{2.7})
and the H\"{o}lder inequality, we may state the following counterparts of
Bessel's inequality as well:%
\begin{multline}
\left\Vert x\right\Vert ^{2}\leq \frac{1}{2}\cdot \frac{1}{\left[ \sum_{i\in
F}\func{Re}\left( \Phi _{i}\overline{\phi _{i}}\right) \right] ^{\frac{1}{2}}%
}  \label{2.11} \\
\times \left\{ 
\begin{array}{l}
\max\limits_{i\in F}\left\{ \left\vert \Phi _{i}\right\vert +\left\vert \phi
_{i}\right\vert \right\} \sum\limits_{i\in F}\left\vert \left\langle
x,e_{i}\right\rangle \right\vert  \\ 
\\ 
\left[ \sum\limits_{i\in F}\left( \left\vert \Phi _{i}\right\vert
+\left\vert \phi _{i}\right\vert \right) ^{p}\right] ^{\frac{1}{p}}\left(
\sum\limits_{i\in F}\left\vert \left\langle x,e_{i}\right\rangle \right\vert
^{q}\right) ^{\frac{1}{q}},\text{ } \\ 
\hfill \text{for\ }p>1,\ \frac{1}{p}+\frac{1}{q}=1 \\ 
\\ 
\max_{i\in F}\left\vert \left\langle x,e_{i}\right\rangle \right\vert
\sum\limits_{i\in F}\left[ \left\vert \Phi _{i}\right\vert +\left\vert \phi
_{i}\right\vert \right] .%
\end{array}%
\right. 
\end{multline}
\end{remark}

The following corollary holds.

\begin{corollary}
\label{c2.3}With the assumption of Theorem \ref{t2.1} and if either (i) or
(ii) holds, then%
\begin{equation}
0\leq \left\Vert x\right\Vert ^{2}-\sum_{i\in F}\left\vert \left\langle
x,e_{i}\right\rangle \right\vert ^{2}\leq \frac{1}{4}M^{2}\left( \mathbf{%
\Phi },\mathbf{\phi },F\right) \sum_{i\in F}\left\vert \left\langle
x,e_{i}\right\rangle \right\vert ^{2},  \label{2.12}
\end{equation}%
where%
\begin{equation}
M\left( \mathbf{\Phi },\mathbf{\phi },F\right) :=\left[ \frac{\sum_{i\in
F}\left\{ \left( \left\vert \Phi _{i}\right\vert +\left\vert \phi
_{i}\right\vert \right) ^{2}+4\left[ \left\vert \Phi _{i}\overline{\phi _{i}}%
\right\vert -\func{Re}\left( \Phi _{i}\overline{\phi _{i}}\right) \right]
\right\} }{\sum_{i\in F}\func{Re}\left( \Phi _{i}\overline{\phi _{i}}\right) 
}\right] ^{\frac{1}{2}}.  \label{2.12.1}
\end{equation}%
The constant $\frac{1}{4}$ is best possible.
\end{corollary}

\begin{proof}
The inequality (\ref{2.12}) follows by (\ref{2.1}) on subtracting the same
quantity $\sum_{i\in F}\left\vert \left\langle x,e_{i}\right\rangle
\right\vert ^{2}$ from both sides.

To prove the sharpness of the constant $\frac{1}{4},$ assume that (\ref{2.12}%
) holds with $c>0,$ i.e., 
\begin{equation}
0\leq \left\Vert x\right\Vert ^{2}-\sum_{i\in F}\left\vert \left\langle
x,e_{i}\right\rangle \right\vert ^{2}\leq cM^{2}\left( \mathbf{\Phi },%
\mathbf{\phi },F\right) \sum_{i\in F}\left\vert \left\langle
x,e_{i}\right\rangle \right\vert ^{2}  \label{2.13}
\end{equation}%
provided the condition (i) holds.

Choose $F=\left\{ 1\right\} ,$ $e_{1}=e,$ $\left\Vert e\right\Vert =1,$ $%
\phi _{i}=\phi ,$ $\Phi _{i}=\Phi ,$ $\phi ,\Phi >0$ in (\ref{2.13}) to get%
\begin{equation}
0\leq \left\Vert x\right\Vert ^{2}-\left\vert \left\langle x,e\right\rangle
\right\vert ^{2}\leq c\frac{\left( \Phi -\phi \right) ^{2}}{\phi \Phi }%
\left\vert \left\langle x,e\right\rangle \right\vert ^{2},  \label{2.14}
\end{equation}%
provided%
\begin{equation}
\left\langle \Phi e-x,x-\phi e\right\rangle \geq 0.  \label{2.14.1}
\end{equation}%
If $H=\mathbb{R}^{2},$ $x=\left( x_{1},x_{2}\right) \in \mathbb{R}^{2},$ $%
e=\left( \frac{1}{\sqrt{2}},\frac{1}{\sqrt{2}}\right) $ then we have%
\begin{align*}
\left\Vert x\right\Vert ^{2}-\left\vert \left\langle x,e\right\rangle
\right\vert ^{2}& =x_{1}^{2}+x_{2}^{2}-\frac{\left( x_{1}+x_{2}\right) ^{2}}{%
2}=\frac{1}{2}\left( x_{1}-x_{2}\right) ^{2}, \\
\left\vert \left\langle x,e\right\rangle \right\vert ^{2}& =\frac{\left(
x_{1}+x_{2}\right) ^{2}}{2}
\end{align*}%
and by (\ref{2.14}) we get%
\begin{equation}
\frac{\left( x_{1}-x_{2}\right) ^{2}}{2}\leq c\frac{\left( \Phi -\phi
\right) ^{2}}{\phi \Phi }\cdot \frac{\left( x_{1}+x_{2}\right) ^{2}}{2}.
\label{2.15}
\end{equation}%
Now, if we let $x_{1}=\frac{\phi }{\sqrt{2}},$ $x_{2}=\frac{\Phi }{\sqrt{2}}$
$\left( \phi ,\Phi >0\right) $ then obviously%
\begin{equation*}
\left\langle \Phi e-x,x-\phi e\right\rangle =\sum_{i=1}^{2}\left( \frac{\Phi 
}{\sqrt{2}}-x_{i}\right) \left( x_{i}-\frac{\phi }{\sqrt{2}}\right) =0,
\end{equation*}%
which shows that (\ref{2.14.1}) is fulfilled, and thus by (\ref{2.15}) we
obtain%
\begin{equation*}
\frac{\left( \Phi -\phi \right) ^{2}}{4}\leq c\frac{\left( \Phi -\phi
\right) ^{2}}{\phi \Phi }\cdot \frac{\left( \Phi +\phi \right) ^{2}}{4}
\end{equation*}%
for any $\Phi >\phi >0.$ This implies%
\begin{equation}
c\left( \Phi +\phi \right) ^{2}\geq \phi \Phi  \label{2.16}
\end{equation}%
for any $\Phi >\phi >0.$

Finally, let $\phi =1-\varepsilon ,$ $\Phi =1+\varepsilon $, $\varepsilon
\in \left( 0,1\right) $. Then from (\ref{2.16}) we get $4c\geq 1-\varepsilon
^{2}$ for any $\varepsilon \in \left( 0,1\right) $ which produces $c\geq 
\frac{1}{4}.$
\end{proof}

\begin{remark}
\label{r2.4}If $\left\{ e_{i}\right\} _{i\in I}$ is an orthornormal family
in the real inner product $\left( H;\left\langle \cdot ,\cdot \right\rangle
\right) $ and $M_{i},m_{i}\in \mathbb{R}$, $i\in F$ ($F$ is a finite part of 
$I$) and $x\in H$ are such that $M_{i},m_{i}\geq 0$ for $i\in F$ with $%
\sum_{i\in F}M_{i}m_{i}\geq 0$ and%
\begin{equation*}
\left\langle \sum_{i\in F}M_{i}e_{i}-x,x-\sum_{i\in
F}m_{i}e_{i}\right\rangle \geq 0,
\end{equation*}%
then we have the inequality%
\begin{equation}
0\leq \left\Vert x\right\Vert ^{2}-\sum_{i\in F}\left[ \left\langle
x,e_{i}\right\rangle \right] ^{2}\leq \frac{1}{4}\cdot \frac{\sum_{i\in
F}\left( M_{i}-m_{i}\right) ^{2}}{\sum_{i\in F}M_{i}m_{i}}\cdot \sum_{i\in F}%
\left[ \left\langle x,e_{i}\right\rangle \right] ^{2}.  \label{2.17}
\end{equation}%
The constant $\frac{1}{4}$ is best possible.
\end{remark}

The following counterpart of the Schwarz's inequality in inner product
spaces holds.

\begin{corollary}
\label{c2.5}Let $x,y\in H$ and $\delta ,\Delta \in \mathbb{K}$ $\left( 
\mathbb{K}=\mathbb{C},\mathbb{R}\right) $ with the property that $\func{Re}%
\left( \Delta \overline{\delta }\right) >0.$ If either%
\begin{equation}
\func{Re}\left\langle \Delta y-x,x-\delta y\right\rangle \geq 0  \label{2.18}
\end{equation}%
or, equivalently,%
\begin{equation}
\left\Vert x-\frac{\delta +\Delta }{2}\cdot y\right\Vert \leq \frac{1}{2}%
\left\vert \Delta -\delta \right\vert \left\Vert y\right\Vert   \label{2.19}
\end{equation}%
holds, then we have the inequalities%
\begin{align}
\left\Vert x\right\Vert \left\Vert y\right\Vert & \leq \frac{1}{2}\cdot 
\frac{\func{Re}\left[ \Delta \overline{\left\langle x,y\right\rangle }+%
\overline{\delta }\left\langle x,y\right\rangle \right] }{\sqrt{\Delta 
\overline{\delta }}}  \label{2.20} \\
& \leq \frac{1}{2}\cdot \frac{\left\vert \Delta \right\vert +\left\vert
\delta \right\vert }{\sqrt{\Delta \overline{\delta }}}\left\vert
\left\langle x,y\right\rangle \right\vert ,  \notag
\end{align}%
\begin{align}
0& \leq \left\Vert x\right\Vert \left\Vert y\right\Vert -\left\vert
\left\langle x,y\right\rangle \right\vert   \label{2.21} \\
& \leq \frac{1}{2}\cdot \frac{\left( \sqrt{\left\vert \Delta \right\vert }-%
\sqrt{\left\vert \delta \right\vert }\right) ^{2}+2\left( \sqrt{\Delta 
\overline{\delta }}-\sqrt{\func{Re}\left( \Delta \overline{\delta }\right) }%
\right) }{\sqrt{\Delta \overline{\delta }}}\left\vert \left\langle
x,y\right\rangle \right\vert ,  \notag
\end{align}%
\begin{equation}
\left\Vert x\right\Vert ^{2}\left\Vert y\right\Vert ^{2}\leq \frac{1}{4}%
\cdot \frac{\left( \left\vert \Delta \right\vert +\left\vert \delta
\right\vert \right) ^{2}}{\func{Re}\left( \Delta \overline{\delta }\right) }%
\left\vert \left\langle x,y\right\rangle \right\vert ^{2},  \label{2.22}
\end{equation}%
and 
\begin{equation}
0\leq \left\Vert x\right\Vert ^{2}\left\Vert y\right\Vert ^{2}-\left\vert
\left\langle x,y\right\rangle \right\vert ^{2}\leq \frac{1}{4}\cdot \frac{%
\left( \left\vert \Delta \right\vert +\left\vert \delta \right\vert \right)
^{2}+4\left( \left\vert \Delta \overline{\delta }\right\vert -\func{Re}%
\left( \Delta \overline{\delta }\right) \right) }{\func{Re}\left( \Delta 
\overline{\delta }\right) }\left\vert \left\langle x,y\right\rangle
\right\vert ^{2}.  \label{2.23}
\end{equation}%
The constants $\frac{1}{2}$ and $\frac{1}{4}$ are best possible.
\end{corollary}

\begin{proof}
The inequality (\ref{2.20}) follows from (\ref{2.6}) on choosing $F=\left\{
1\right\} ,$ $e_{1}=e=\frac{y}{\left\Vert y\right\Vert },$ $\Phi _{1}=\Phi
=\Delta \left\Vert y\right\Vert ,$ \ $\phi _{1}=\phi =\delta \left\Vert
y\right\Vert $ $\left( y\neq 0\right) .$ The inequality (\ref{2.21}) is
equivalent with (\ref{2.20}). The inequality (\ref{2.22}) follows from (\ref%
{2.1}) for $F=\left\{ 1\right\} $ and the same choices as above. Finally, (%
\ref{2.23}) is obviously equivalent with (\ref{2.22}).
\end{proof}

\section{Some Gr\"{u}ss Type Inequalities\label{s3}}

The following result holds.

\begin{theorem}
\label{t3.1}Let $\left\{ e_{i}\right\} _{i\in I}$ be a family of
orthornormal vectors in $H,$ $F$ a finite part of $I$, $\phi _{i},\Phi _{i},$
$\gamma _{i},\Gamma _{i}\in \mathbb{K},\ i\in F$ and $x,y\in H.$ If either%
\begin{align}
\func{Re}\left\langle \sum_{i\in F}\Phi _{i}e_{i}-x,x-\sum_{i\in F}\phi
_{i}e_{i}\right\rangle & \geq 0,  \label{3.1} \\
\func{Re}\left\langle \sum_{i\in F}\Gamma _{i}e_{i}-y,y-\sum_{i\in F}\gamma
_{i}e_{i}\right\rangle & \geq 0,  \notag
\end{align}%
or, equivalently,%
\begin{align}
\left\Vert x-\sum_{i\in F}\frac{\Phi _{i}+\phi _{i}}{2}e_{i}\right\Vert &
\leq \frac{1}{2}\left( \sum_{i\in F}\left\vert \Phi _{i}-\phi
_{i}\right\vert ^{2}\right) ^{\frac{1}{2}},  \label{3.2} \\
\left\Vert y-\sum_{i\in F}\frac{\Gamma _{i}+\gamma _{i}}{2}e_{i}\right\Vert
& \leq \frac{1}{2}\left( \sum_{i\in F}\left\vert \Gamma _{i}-\gamma
_{i}\right\vert ^{2}\right) ^{\frac{1}{2}},  \notag
\end{align}%
hold, then we have the inequality%
\begin{align}
0& \leq \left\vert \left\langle x,y\right\rangle -\sum_{i\in F}\left\langle
x,e_{i}\right\rangle \left\langle e_{i},y\right\rangle \right\vert
\label{3.3} \\
& \leq \frac{1}{4}M\left( \mathbf{\Phi },\mathbf{\phi },F\right) M\left( 
\mathbf{\Gamma },\mathbf{\gamma },F\right) \left( \sum_{i\in F}\left\vert
\left\langle x,e_{i}\right\rangle \right\vert ^{2}\right) ^{\frac{1}{2}%
}\left( \sum_{i\in F}\left\vert \left\langle y,e_{i}\right\rangle
\right\vert ^{2}\right) ^{\frac{1}{2}},  \notag
\end{align}%
where $M\left( \mathbf{\Phi },\mathbf{\phi },F\right) $ is defined in (\ref%
{2.12.1}).

The constant $\frac{1}{4}$ is best possible.
\end{theorem}

\begin{proof}
Using Schwartz's inequality in the inner product space $\left(
H,\left\langle \cdot ,\cdot \right\rangle \right) $ one has%
\begin{multline}
\left\vert \left\langle x-\sum_{i\in F}\left\langle x,e_{i}\right\rangle
e_{i},y-\sum_{i\in F}\left\langle y,e_{i}\right\rangle e_{i}\right\rangle
\right\vert ^{2}  \label{3.4} \\
\leq \left\Vert x-\sum_{i\in F}\left\langle x,e_{i}\right\rangle
e_{i}\right\Vert ^{2}\left\Vert y-\sum_{i\in F}\left\langle
y,e_{i}\right\rangle e_{i}\right\Vert ^{2}
\end{multline}%
and since a simple calculation shows that 
\begin{equation*}
\left\langle x-\sum_{i\in F}\left\langle x,e_{i}\right\rangle
e_{i},y-\sum_{i\in F}\left\langle y,e_{i}\right\rangle e_{i}\right\rangle
=\left\langle x,y\right\rangle -\sum_{i\in F}\left\langle
x,e_{i}\right\rangle \left\langle e_{i},y\right\rangle
\end{equation*}%
and 
\begin{equation*}
\left\Vert x-\sum_{i\in F}\left\langle x,e_{i}\right\rangle e_{i}\right\Vert
^{2}\leq \left\Vert x\right\Vert ^{2}-\sum_{i\in F}\left\vert \left\langle
x,e_{i}\right\rangle \right\vert ^{2}
\end{equation*}%
for any $x,y\in H,$ then by (\ref{3.4}) and by the counterpart of Bessel's
inequality in Corollary \ref{c2.3}, we have%
\begin{align}
& \left\vert \left\langle x,y\right\rangle -\sum_{i\in F}\left\langle
x,e_{i}\right\rangle \left\langle e_{i},y\right\rangle \right\vert ^{2}
\label{3.5} \\
& \leq \left( \left\Vert x\right\Vert ^{2}-\sum_{i\in F}\left\vert
\left\langle x,e_{i}\right\rangle \right\vert ^{2}\right) \left( \left\Vert
y\right\Vert ^{2}-\sum_{i\in F}\left\vert \left\langle y,e_{i}\right\rangle
\right\vert ^{2}\right)  \notag \\
& \leq \frac{1}{4}M^{2}\left( \mathbf{\Phi },\mathbf{\phi },F\right)
\sum_{i\in F}\left\vert \left\langle x,e_{i}\right\rangle \right\vert
^{2}\cdot \frac{1}{4}M^{2}\left( \mathbf{\Gamma },\mathbf{\gamma },F\right)
\sum_{i\in F}\left\vert \left\langle y,e_{i}\right\rangle \right\vert ^{2}. 
\notag
\end{align}%
Taking the square root in (\ref{3.5}), we deduce (\ref{3.3}).

The fact that $\frac{1}{4}$ is the best possible constant follows by
Corollary \ref{c2.3} and we omit the details.
\end{proof}

The following corollary for real inner product spaces holds.

\begin{corollary}
\label{c3.2}Let $\left\{ e_{i}\right\} _{i\in I}$ be a family of
orthornormal vectors in $H,$ $F$ a finite part of $I$, $M_{i},m_{i},$ $%
N_{i},n_{i}\geq 0,\ i\in F$ and $x,y\in H$ such that $\sum_{i\in
F}M_{i}m_{i}>0,$ $\sum_{i\in F}N_{i}n_{i}>0$ and%
\begin{equation}
\left\langle \sum_{i\in F}M_{i}e_{i}-x,x-\sum_{i\in
F}m_{i}e_{i}\right\rangle \geq 0,\ \ \ \ \left\langle \sum_{i\in
F}N_{i}e_{i}-y,y-\sum_{i\in F}n_{i}e_{i}\right\rangle \geq 0.  \label{3.6}
\end{equation}%
Then we have the inequality%
\begin{align}
0& \leq \left\vert \left\langle x,y\right\rangle -\sum_{i\in F}\left\langle
x,e_{i}\right\rangle \left\langle y,e_{i}\right\rangle \right\vert ^{2}
\label{3.7} \\
& \leq \frac{1}{16}\cdot \frac{\sum_{i\in F}\left( M_{i}-m_{i}\right)
^{2}\sum_{i\in F}\left( N_{i}-n_{i}\right) ^{2}\sum_{i\in F}\left\vert
\left\langle x,e_{i}\right\rangle \right\vert ^{2}\sum_{i\in F}\left\vert
\left\langle y,e_{i}\right\rangle \right\vert ^{2}}{\sum_{i\in
F}M_{i}m_{i}\sum_{i\in F}N_{i}n_{i}}.  \notag
\end{align}%
The constant $\frac{1}{16}$ is best possible.
\end{corollary}

In the case where the family $\left\{ e_{i}\right\} _{i\in I}$ reduces to a
single vector, we may deduce from Theorem \ref{t3.1} the following
particular case first obtained in \cite{SSDc}.

\begin{corollary}
\label{c3.3}Let $e\in H,$ $\left\Vert e\right\Vert =1,$ $\phi ,\Phi ,\gamma
,\Gamma \in \mathbb{K}$ \ with $\func{Re}\left( \Phi \overline{\phi }\right)
,$ $\func{Re}\left( \Gamma \overline{\gamma }\right) >0$ and $x,y\in H$ such
that either 
\begin{equation}
\func{Re}\left\langle \Phi e-x,x-\phi e\right\rangle \geq 0,\ \ \func{Re}%
\left\langle \Gamma e-y,y-\gamma e\right\rangle \geq 0,  \label{3.8}
\end{equation}%
or, equivalently,%
\begin{equation}
\left\Vert x-\frac{\phi +\Phi }{2}e\right\Vert \leq \frac{1}{2}\left\vert
\Phi -\phi \right\vert ,\ \ \ \ \ \left\Vert y-\frac{\gamma +\Gamma }{2}%
e\right\Vert \leq \frac{1}{2}\left\vert \Gamma -\gamma \right\vert 
\label{3.9}
\end{equation}%
holds, then%
\begin{equation}
0\leq \left\vert \left\langle x,y\right\rangle -\left\langle
x,e\right\rangle \left\langle e,y\right\rangle \right\vert \leq \frac{1}{4}%
M\left( \Phi ,\phi \right) M\left( \Gamma ,\gamma \right) \left\vert
\left\langle x,e\right\rangle \left\langle e,y\right\rangle \right\vert ,
\label{3.10}
\end{equation}%
where%
\begin{equation*}
M\left( \Phi ,\phi \right) :=\left[ \frac{\left( \left\vert \Phi \right\vert
-\left\vert \phi \right\vert \right) ^{2}+4\left[ \left\vert \phi \Phi
\right\vert -\func{Re}\left( \Phi \overline{\phi }\right) \right] }{\func{Re}%
\left( \Phi \overline{\phi }\right) }\right] ^{\frac{1}{2}}.
\end{equation*}%
The constant $\frac{1}{4}$ is best possible.
\end{corollary}

\begin{remark}
\label{r3.4}If $H$ is real, $e\in H,$ $\left\Vert e\right\Vert =1$ and $%
a,b,A,B\in \mathbb{R}$ are such that $A>a>0,$ $B>b>0$ and%
\begin{equation}
\left\Vert x-\frac{a+A}{2}e\right\Vert \leq \frac{1}{2}\left( A-a\right) ,\
\ \left\Vert y-\frac{b+B}{2}e\right\Vert \leq \frac{1}{2}\left( B-b\right) ,
\label{3.11}
\end{equation}%
then%
\begin{equation}
\left\vert \left\langle x,y\right\rangle -\left\langle x,e\right\rangle
\left\langle e,y\right\rangle \right\vert \leq \frac{1}{4}\cdot \frac{\left(
A-a\right) \left( B-b\right) }{\sqrt{abAB}}\left\vert \left\langle
x,e\right\rangle \left\langle e,y\right\rangle \right\vert .  \label{3.12}
\end{equation}%
The constant $\frac{1}{4}$ is best possible.
\end{remark}

If $\left\langle x,e\right\rangle ,$ $\left\langle y,e\right\rangle \neq 0,$
then the following equivalent form of (\ref{3.12}) also holds%
\begin{equation}
\left\vert \frac{\left\langle x,y\right\rangle }{\left\langle
x,e\right\rangle \left\langle e,y\right\rangle }-1\right\vert \leq \frac{1}{4%
}\cdot \frac{\left( A-a\right) \left( B-b\right) }{\sqrt{abAB}}.
\label{3.13}
\end{equation}

\section{Some Companion Inequalities\label{s4}}

The following companion of the Gr\"{u}ss inequality also holds.

\begin{theorem}
\label{t4.1}Let $\left\{ e_{i}\right\} _{i\in I}$ be a family of
orthornormal vectors in $H,$ $F$ a finite part of $I$, $\phi _{i},\Phi
_{i}\in \mathbb{K},\ \left( i\in F\right) $, $x,y\in H$ and $\lambda \in
\left( 0,1\right) ,$ such that either%
\begin{equation}
\func{Re}\left\langle \sum_{i\in F}\Phi _{i}e_{i}-\left( \lambda x+\left(
1-\lambda \right) y\right) ,\lambda x+\left( 1-\lambda \right) y-\sum_{i\in
F}\phi _{i}e_{i}\right\rangle \geq 0  \label{4.1}
\end{equation}%
or, equivalently,%
\begin{equation}
\left\Vert \lambda x+\left( 1-\lambda \right) y-\sum_{i\in F}\frac{\Phi
_{i}+\phi _{i}}{2}\cdot e_{i}\right\Vert \leq \frac{1}{2}\left( \sum_{i\in
F}\left\vert \Phi _{i}-\phi _{i}\right\vert ^{2}\right) ^{\frac{1}{2}},
\label{4.2}
\end{equation}%
holds. Then we have the inequality%
\begin{multline}
\func{Re}\left[ \left\langle x,y\right\rangle -\sum_{i\in F}\left\langle
x,e_{i}\right\rangle \left\langle e_{i},y\right\rangle \right]  \label{4.3}
\\
\leq \frac{1}{16}\cdot \frac{1}{\lambda \left( 1-\lambda \right) }\sum_{i\in
F}M^{2}\left( \mathbf{\Phi },\mathbf{\phi },F\right) \sum_{i\in F}\left\vert
\left\langle \lambda x+\left( 1-\lambda \right) y,e_{i}\right\rangle
\right\vert ^{2}.
\end{multline}%
The constant $\frac{1}{16}$ is the best possible constant in (\ref{4.3}) in
the sense that it cannot be replaced by a smaller constant.
\end{theorem}

\begin{proof}
Using the known inequality%
\begin{equation*}
\func{Re}\left\langle z,u\right\rangle \leq \frac{1}{4}\left\Vert
z+u\right\Vert ^{2}
\end{equation*}%
we may state that for any $a,b\in H$ and $\lambda \in \left( 0,1\right) $ 
\begin{equation}
\func{Re}\left\langle a,b\right\rangle \leq \frac{1}{4\lambda \left(
1-\lambda \right) }\left\Vert \lambda a+\left( 1-\lambda \right)
b\right\Vert ^{2}.  \label{4.4}
\end{equation}%
Since%
\begin{equation*}
\left\langle x,y\right\rangle -\sum_{i\in F}\left\langle
x,e_{i}\right\rangle \left\langle e_{i},y\right\rangle =\left\langle
x-\sum_{i\in F}\left\langle x,e_{i}\right\rangle e_{i},y-\sum_{i\in
F}\left\langle y,e_{i}\right\rangle e_{i}\right\rangle ,
\end{equation*}%
for any \thinspace $x,y\in H,$ then, by (\ref{4.4}), we get%
\begin{align}
& \func{Re}\left[ \left\langle x,y\right\rangle -\sum_{i\in F}\left\langle
x,e_{i}\right\rangle \left\langle e_{i},y\right\rangle \right]  \label{4.5}
\\
& =\func{Re}\left[ \left\langle x-\sum_{i\in F}\left\langle
x,e_{i}\right\rangle e_{i},y-\sum_{i\in F}\left\langle y,e_{i}\right\rangle
e_{i}\right\rangle \right]  \notag \\
& \leq \frac{1}{4\lambda \left( 1-\lambda \right) }\left\Vert \lambda \left(
x-\sum_{i\in F}\left\langle x,e_{i}\right\rangle e_{i}\right) +\left(
1-\lambda \right) \left( y-\sum_{i\in F}\left\langle y,e_{i}\right\rangle
e_{i}\right) \right\Vert ^{2}  \notag \\
& =\frac{1}{4\lambda \left( 1-\lambda \right) }\left\Vert \lambda x+\left(
1-\lambda \right) y-\sum_{i\in F}\left\langle \lambda x+\left( 1-\lambda
\right) y,e_{i}\right\rangle e_{i}\right\Vert ^{2}  \notag \\
& =\frac{1}{4\lambda \left( 1-\lambda \right) }\left[ \left\Vert \lambda
x+\left( 1-\lambda \right) y\right\Vert ^{2}-\sum_{i\in F}\left\vert
\left\langle \lambda x+\left( 1-\lambda \right) y,e_{i}\right\rangle
\right\vert ^{2}\right] .  \notag
\end{align}%
If we apply the counterpart of Bessel's inequality from Corollary \ref{c2.3}
for $\lambda x+\left( 1-\lambda \right) y,$ we may state that%
\begin{multline}
\left\Vert \lambda x+\left( 1-\lambda \right) y\right\Vert ^{2}-\sum_{i\in
F}\left\vert \left\langle \lambda x+\left( 1-\lambda \right)
y,e_{i}\right\rangle \right\vert ^{2}  \label{4.6} \\
\leq \frac{1}{4}M^{2}\left( \mathbf{\Phi },\mathbf{\phi },F\right)
\sum_{i\in F}\left\vert \left\langle \lambda x+\left( 1-\lambda \right)
y,e_{i}\right\rangle \right\vert ^{2}.
\end{multline}%
Now, by making use of (\ref{4.5}) and (\ref{4.6}), we deduce (\ref{4.3}).

The fact that $\frac{1}{16}$ is the best possible constant in (\ref{4.3})
follows by the fact that if in (\ref{4.1}) we choose $x=y,$ then it becomes
(i) of Theorem \ref{t2.1}, implying for $\lambda =\frac{1}{2}$ (\ref{2.12}),
for which, we have shown that $\frac{1}{4}$ was the best constant.
\end{proof}

\begin{remark}
\label{r4.2}If in Theorem \ref{t4.1}, we choose $\lambda =\frac{1}{2},$ then
we get%
\begin{equation}
\func{Re}\left[ \left\langle x,y\right\rangle -\sum_{i\in F}\left\langle
x,e_{i}\right\rangle \left\langle e_{i},y\right\rangle \right] \leq \frac{1}{%
4}M^{2}\left( \mathbf{\Phi },\mathbf{\phi },F\right) \sum_{i\in F}\left\vert
\left\langle \frac{x+y}{2},e_{i}\right\rangle \right\vert ^{2},  \label{4.7}
\end{equation}%
provided%
\begin{equation*}
\func{Re}\left\langle \sum_{i\in F}\Phi _{i}e_{i}-\frac{x+y}{2},\frac{x+y}{2}%
-\sum_{i\in F}\phi _{i}e_{i}\right\rangle \geq 0
\end{equation*}%
or, equivalently,%
\begin{equation}
\left\Vert \frac{x+y}{2}-\sum_{i\in F}\frac{\Phi _{i}+\phi _{i}}{2}\cdot
e_{i}\right\Vert \leq \frac{1}{2}\left( \sum_{i\in F}\left\vert \Phi
_{i}-\phi _{i}\right\vert ^{2}\right) ^{\frac{1}{2}}.  \label{4.8}
\end{equation}
\end{remark}

\section{Integral Inequalities\label{s5}}

Let $\left( \Omega ,\Sigma ,\mu \right) $ be a measure space consisting of a
set $\Omega ,$ a $\sigma -$algebra of parts $\Sigma $ and a countably
additive and positive measure $\mu $ on $\Sigma $ with values in $\mathbb{R}%
\cup \left\{ \infty \right\} .$ Let $\rho \geq 0$ be a $\mu -$measurable
function on $\Omega .$ Denote by $L_{\rho }^{2}\left( \Omega ,\mathbb{K}%
\right) $ the Hilbert space of all real or complex valued functions defined
on $\Omega $ and $2-\rho -$integrable on $\Omega ,$ i.e.,%
\begin{equation}
\int_{\Omega }\rho \left( s\right) \left\vert f\left( s\right) \right\vert
^{2}d\mu \left( s\right) <\infty .  \label{5.1}
\end{equation}

Consider the family $\left\{ f_{i}\right\} _{i\in I}$ of functions in $%
L_{\rho }^{2}\left( \Omega ,\mathbb{K}\right) $ with the properties that%
\begin{equation}
\int_{\Omega }\rho \left( s\right) f_{i}\left( s\right) \overline{f_{j}}%
\left( s\right) d\mu \left( s\right) =\delta _{ij},\ \ \ i,j\in I,
\label{5.2}
\end{equation}%
where $\delta _{ij}$ is $0$ if $i\neq j$ and $\delta _{ij}=1$ if $i=j.$ $%
\left\{ f_{i}\right\} _{i\in I}$ is an orthornormal family in $L_{\rho
}^{2}\left( \Omega ,\mathbb{K}\right) .$

The following proposition holds.

\begin{proposition}
\label{p5.1}Let $\left\{ f_{i}\right\} _{i\in I}$ be an orthornormal family
of functions in $L_{\rho }^{2}\left( \Omega ,\mathbb{K}\right) ,$ $F$ a
finite subset of $I,$ $\phi _{i},\Phi _{i}\in \mathbb{K}$ $\left( i\in
F\right) $ such that $\sum_{i\in F}\func{Re}\left( \Phi _{i}\overline{\phi
_{i}}\right) >0$ and $f\in L_{\rho }^{2}\left( \Omega ,\mathbb{K}\right) ,$
so that either%
\begin{equation}
\int_{\Omega }\rho \left( s\right) \func{Re}\left[ \left( \sum_{i\in F}\Phi
_{i}f_{i}\left( s\right) -f\left( s\right) \right) \left( \overline{f}\left(
s\right) -\sum_{i\in F}\overline{\phi _{i}}\text{ }\overline{f_{i}}\left(
s\right) \right) \right] d\mu \left( s\right) \geq 0  \label{5.3}
\end{equation}%
or, equivalently,%
\begin{equation}
\int_{\Omega }\rho \left( s\right) \left\vert f\left( s\right) -\sum_{i\in F}%
\frac{\Phi _{i}+\phi _{i}}{2}f_{i}\left( s\right) \right\vert ^{2}d\mu
\left( s\right) \leq \frac{1}{4}\sum_{i\in F}\left\vert \Phi _{i}-\phi
_{i}\right\vert ^{2}.  \label{5.4}
\end{equation}%
Then we have the inequality%
\begin{equation}
\left( \int_{\Omega }\rho \left( s\right) \left\vert f\left( s\right)
\right\vert ^{2}d\mu \left( s\right) \right) ^{\frac{1}{2}}\leq \frac{1}{2}%
\cdot \frac{1}{\left[ \sum_{i\in F}\func{Re}\left( \Phi _{i}\overline{\phi
_{i}}\right) \right] ^{\frac{1}{2}}}  \label{5.5}
\end{equation}%
\begin{equation*}
\times \left\{ 
\begin{array}{l}
\max\limits_{i\in F}\left\{ \left\vert \Phi _{i}\right\vert +\left\vert \phi
_{i}\right\vert \right\} \dsum\limits_{i\in F}\left\vert \dint_{\Omega }\rho
\left( s\right) f\left( s\right) \overline{f_{i}}\left( s\right) d\mu \left(
s\right) \right\vert  \\ 
\\ 
\left[ \dsum\limits_{i\in F}\left( \left\vert \Phi _{i}\right\vert
+\left\vert \phi _{i}\right\vert \right) ^{p}\right] ^{\frac{1}{p}}\left(
\dsum\limits_{i\in F}\left\vert \dint_{\Omega }\rho \left( s\right) f\left(
s\right) \overline{f_{i}}\left( s\right) d\mu \left( s\right) \right\vert
^{q}\right) ^{\frac{1}{q}},\text{ } \\ 
\hfill \text{\ for \ }p>1,\ \frac{1}{p}+\frac{1}{q}=1 \\ 
\\ 
\max\limits_{i\in F}\left\vert \dint_{\Omega }\rho \left( s\right) f\left(
s\right) \overline{f_{i}}\left( s\right) d\mu \left( s\right) \right\vert
\dsum\limits_{i\in F}\left[ \left\vert \Phi _{i}\right\vert +\left\vert \phi
_{i}\right\vert \right] .%
\end{array}%
\right. 
\end{equation*}%
In particular, we have%
\begin{multline}
\int_{\Omega }\rho \left( s\right) \left\vert f\left( s\right) \right\vert
^{2}d\mu \left( s\right)   \label{5.6} \\
\leq \frac{1}{4}\cdot \frac{\sum_{i\in F}\left( \left\vert \Phi
_{i}\right\vert +\left\vert \phi _{i}\right\vert \right) ^{2}}{\sum_{i\in F}%
\func{Re}\left( \Phi _{i}\overline{\phi _{i}}\right) }\sum\limits_{i\in
F}\left\vert \int_{\Omega }\rho \left( s\right) f\left( s\right) \overline{%
f_{i}}\left( s\right) d\mu \left( s\right) \right\vert ^{2}.
\end{multline}%
The constant $\frac{1}{4}$ is best possible in both inequalities.
\end{proposition}

The proof is obvious by Theorem \ref{t2.1} and Remark \ref{r2.2}. We omit
the details.

The following proposition also holds.

\begin{proposition}
\label{p5.2}Assume that $f_{i},f,\phi _{i},\Phi _{i}$ and $F$ satisfy the
assumptions of Proposition \ref{p5.1}. Then we have the following
counterpart of Bessel's inequality:%
\begin{align}
0& \leq \int_{\Omega }\rho \left( s\right) f^{2}\left( s\right) d\mu \left(
s\right) -\sum\limits_{i\in F}\left\vert \int_{\Omega }\rho \left( s\right)
f\left( s\right) \overline{f_{i}}\left( s\right) d\mu \left( s\right)
\right\vert ^{2}  \label{5.7} \\
& \leq \frac{1}{4}M^{2}\left( \mathbf{\Phi },\mathbf{\phi },F\right) \cdot
\sum\limits_{i\in F}\left\vert \int_{\Omega }\rho \left( s\right) f\left(
s\right) \overline{f_{i}}\left( s\right) d\mu \left( s\right) \right\vert
^{2},  \notag
\end{align}%
where, as above,%
\begin{equation}
M\left( \mathbf{\Phi },\mathbf{\phi },F\right) :=\left[ \frac{%
\sum\limits_{i\in F}\left\{ \left( \left\vert \Phi _{i}\right\vert
-\left\vert \phi _{i}\right\vert \right) ^{2}+4\left[ \left\vert \phi
_{i}\Phi _{i}\right\vert -\func{Re}\left( \Phi _{i}\overline{\phi _{i}}%
\right) \right] \right\} }{\func{Re}\left( \Phi _{i}\overline{\phi _{i}}%
\right) }\right] ^{\frac{1}{2}}.  \label{5.8}
\end{equation}%
The constant $\frac{1}{4}$ is the best possible.
\end{proposition}

The following Gr\"{u}ss type inequality also holds.

\begin{proposition}
\label{p5.3}Let $\left\{ f_{i}\right\} _{i\in I}$ and $F$ be as in
Proposition \ref{p5.1}. If $\phi _{i},\Phi _{i},\gamma _{i},\Gamma _{i}\in 
\mathbb{K}$ $\left( i\in F\right) $ and $f,g\in L_{\rho }^{2}\left( \Omega ,%
\mathbb{K}\right) $ so that either%
\begin{align}
\int_{\Omega }\rho \left( s\right) \func{Re}\left[ \left( \sum_{i\in F}\Phi
_{i}f_{i}\left( s\right) -f\left( s\right) \right) \left( \overline{f}\left(
s\right) -\sum_{i\in F}\overline{\phi _{i}}\text{ }\overline{f_{i}}\left(
s\right) \right) \right] d\mu \left( s\right) & \geq 0,  \label{5.9} \\
\int_{\Omega }\rho \left( s\right) \func{Re}\left[ \left( \sum_{i\in
F}\Gamma _{i}f_{i}\left( s\right) -g\left( s\right) \right) \left( \overline{%
g}\left( s\right) -\sum_{i\in F}\overline{\gamma _{i}}\text{ }\overline{f_{i}%
}\left( s\right) \right) \right] d\mu \left( s\right) & \geq 0,  \notag
\end{align}%
or, equivalently,%
\begin{align}
\int_{\Omega }\rho \left( s\right) \left\vert f\left( s\right) -\sum_{i\in F}%
\frac{\Phi _{i}+\phi _{i}}{2}\cdot f_{i}\left( s\right) \right\vert ^{2}d\mu
\left( s\right) & \leq \frac{1}{4}\sum_{i\in F}\left\vert \Phi _{i}-\phi
_{i}\right\vert ^{2},  \label{5.10} \\
\int_{\Omega }\rho \left( s\right) \left\vert g\left( s\right) -\sum_{i\in F}%
\frac{\Gamma _{i}+\gamma _{i}}{2}\cdot f_{i}\left( s\right) \right\vert
^{2}d\mu \left( s\right) & \leq \frac{1}{4}\sum_{i\in F}\left\vert \Gamma
_{i}-\gamma _{i}\right\vert ^{2},  \notag
\end{align}%
then we have the inequality%
\begin{multline}
\left\vert \int_{\Omega }\rho \left( s\right) f\left( s\right) \overline{%
g\left( s\right) }d\mu \left( s\right) \right.   \label{5.11} \\
-\left. \sum_{i\in F}\int_{\Omega }\rho \left( s\right) f\left( s\right) 
\overline{f_{i}}\left( s\right) d\mu \left( s\right) \int_{\Omega }\rho
\left( s\right) f_{i}\left( s\right) \overline{g\left( s\right) }d\mu \left(
s\right) \right\vert  \\
\leq \frac{1}{4}M\left( \mathbf{\Phi },\mathbf{\phi },F\right) M\left( 
\mathbf{\Gamma },\mathbf{\gamma },F\right) \left( \sum_{i\in F}\left\vert
\int_{\Omega }\rho \left( s\right) f\left( s\right) \overline{f_{i}}\left(
s\right) d\mu \left( s\right) \right\vert ^{2}\right) ^{\frac{1}{2}} \\
\times \left( \sum_{i\in F}\left\vert \rho \left( s\right) f_{i}\left(
s\right) \overline{g\left( s\right) }d\mu \left( s\right) \right\vert
^{2}\right) ^{\frac{1}{2}},
\end{multline}%
where $M\left( \mathbf{\Phi },\mathbf{\phi },F\right) $ is defined in (\ref%
{5.8}).

The constant $\frac{1}{4}$ is the best possible.
\end{proposition}

The proof follows by Theorem \ref{t3.1} and we omit the details.

In the case of real spaces, the following corollaries provide much simpler
sufficient conditions for the counterpart of Bessel's inequality (\ref{5.7})
or for the Gr\"{u}ss type inequality (\ref{5.11}) to hold.

\begin{corollary}
\label{c5.4}Let $\left\{ f_{i}\right\} _{i\in I}$ be an orthornormal family
of functions in the real Hilbert space $L_{\rho }^{2}\left( \Omega \right) ,$
$F$ a finite part of $I,$ $M_{i},m_{i}\geq 0$ \ $\left( i\in F\right) ,$
with $\sum_{i\in F}M_{i}m_{i}>0$ and $f\in L_{\rho }^{2}\left( \Omega
\right) $ so that%
\begin{equation}
\sum_{i\in F}m_{i}f_{i}\left( s\right) \leq f\left( s\right) \leq \sum_{i\in
F}M_{i}f_{i}\left( s\right) \text{ \ for \ }\mu -\text{a.e. }s\in \Omega .
\label{5.12}
\end{equation}%
Then we have the inequality%
\begin{align}
0& \leq \int_{\Omega }\rho \left( s\right) f^{2}\left( s\right) d\mu \left(
s\right) -\sum_{i\in F}\left[ \int_{\Omega }\rho \left( s\right) f\left(
s\right) f_{i}\left( s\right) d\mu \left( s\right) \right] ^{2}  \label{5.13}
\\
& \leq \frac{1}{4}\cdot \frac{\sum_{i\in F}\left( M_{i}-m_{i}\right) ^{2}}{%
\sum_{i\in F}M_{i}m_{i}}\cdot \sum_{i\in F}\left[ \int_{\Omega }\rho \left(
s\right) f\left( s\right) f_{i}\left( s\right) d\mu \left( s\right) \right]
^{2}.  \notag
\end{align}%
The constant $\frac{1}{4}$ is best possible.
\end{corollary}

\begin{corollary}
\label{c5.5}Let $\left\{ f_{i}\right\} _{i\in I}$ and $F$ be as above. If $%
M_{i},m_{i},N_{i},n_{i}\geq 0$ $\left( i\in F\right) $ with $\sum_{i\in
F}M_{i}m_{i},\sum_{i\in F}N_{i}n_{i}>0$ and $f,g\in L_{\rho }^{2}\left(
\Omega \right) $ such that%
\begin{equation}
\sum_{i\in F}m_{i}f_{i}\left( s\right) \leq f\left( s\right) \leq \sum_{i\in
F}M_{i}f_{i}\left( s\right)   \label{5.14}
\end{equation}%
and%
\begin{equation*}
\sum_{i\in F}n_{i}f_{i}\left( s\right) \leq g\left( s\right) \leq \sum_{i\in
F}N_{i}f_{i}\left( s\right) \text{ \ for \ }\mu -\text{a.e. }s\in \Omega ,
\end{equation*}%
then we have the inequality%
\begin{multline}
\left\vert \int_{\Omega }\rho \left( s\right) f\left( s\right) g\left(
s\right) d\mu \left( s\right) \right.   \label{5.15} \\
-\left. \sum_{i\in F}\int_{\Omega }\rho \left( s\right) f\left( s\right)
f_{i}\left( s\right) d\mu \left( s\right) \int_{\Omega }\rho \left( s\right)
g\left( s\right) f_{i}\left( s\right) d\mu \left( s\right) \right\vert  \\
\leq \frac{1}{4}\left( \frac{\sum_{i\in F}\left( M_{i}-m_{i}\right) ^{2}}{%
\sum_{i\in F}M_{i}m_{i}}\right) ^{\frac{1}{2}}\left( \frac{\sum_{i\in
F}\left( N_{i}-n_{i}\right) ^{2}}{\sum_{i\in F}N_{i}n_{i}}\right) ^{\frac{1}{%
2}} \\
\times \left[ \sum_{i\in F}\left( \int_{\Omega }\rho \left( s\right) f\left(
s\right) f_{i}\left( s\right) d\mu \left( s\right) \right) ^{2}\sum_{i\in
F}\left( \int_{\Omega }\rho \left( s\right) g\left( s\right) f_{i}\left(
s\right) d\mu \left( s\right) \right) ^{2}\right] ^{\frac{1}{2}}.  \notag
\end{multline}%
\qquad 
\end{corollary}


\begin{thebibliography}{9}
\bibitem{SSDa} S.S. DRAGOMIR, A counterpart of Bessel's inequality in inner
product spaces and some Gr\"{u}ss type related results, \textit{RGMIA Res.
Rep. Coll.}, \textbf{6}(2003), Supplement, Article 10. \texttt{[ON Line:
http://rgmia.vu.edu.au/v6(E).html]}

\bibitem{SSDb} S.S. DRAGOMIR, On Bessel and Gr\"{u}ss inequalities for
orthornormal families in inner product spaces, \textit{RGMIA Res. Rep. Coll.}%
, \textbf{6}(2003), Supplement, \texttt{[ON Line:
http://rgmia.vu.edu.au/v6(E).html]}

\bibitem{SSDc} S.S. DRAGOMIR, Some Gr\"{u}ss' type inequalities in inner
product spaces, \textit{J. Ineq. Pure \& Appl. Math., }to appear. [ONLINE: 
\texttt{http://jipam.vu.edu.au}]
\end{thebibliography}
\end{document}